\documentstyle[amssymb,12pt]{amsart}

\newtheorem{Theorem}{Theorem}[section]
\newtheorem{Question}[Theorem]{Question}
\newtheorem{Lemma}[Theorem]{Lemma}

\newtheorem{Remark}[Theorem]{Remark}

\newcommand {\bom}{{\omega^{\Bbb B}}}
\newcommand {\boul}{{\omega^{\Bbb B}/{\cal E}}}
\newcommand {\ul}{{\omega^{\kappa}/{\cal F}}}
\newcommand {\boB}{{\Bbb B}}
\newcommand {\boC}{{\Bbb C}}
\newcommand {\kl}{{{\cal P}_{\kappa}(\lambda)}}
\newcommand {\F}{{\cal F}}
\newcommand {\D}{{\cal D}}
\newcommand {\E}{{\cal E}}
\newcommand {\C}{{\cal C}}

\newcommand{\res}{\upharpoonright}

\begin{document}

\baselineskip=24pt

  \begin{center}
     {\large\bf Possible Size of an Ultrapower of $\omega$} 
  \end{center}

  \begin{center}
 Renling Jin\footnote{The research of the first author
was supported by NSF postdoctoral fellowship \#DMS-9508887.} 
\& Saharon Shelah\footnote{The research
of the second author was supported by The Israel Science Foundation
administered by The Israel Academy of Sciences and Humanities. 
This paper is number 626 on the second author's publication list.}
  \end{center}

  \bigskip

  \begin{quote}

    \centerline{Abstract}

    \small  
Let $\omega$ be the first infinite ordinal (or the set of all natural
numbers) with the usual order $<$. In \S 1 we show that, assuming the 
consistency of a supercompact cardinal, there may exist an ultrapower
of $\omega$, whose cardinality is (1) a singular strong limit cardinal,
(2) a strongly inaccessible cardinal. This answers two questions
in \cite{CK}, modulo the assumption of supercompactness. 
In \S 2 we construct several 
$\lambda$-Archimedean ultrapowers of $\omega$ under some large cardinal 
assumptions. For example, we show that, assuming the 
consistency of a measurable cardinal, there may exist a 
$\lambda$-Archimedean ultrapower of $\omega$ for some 
uncountable cardinal $\lambda$. This answers a question in \cite{KS}, modulo 
the assumption of measurability.
 
\end{quote}

\section{On Notation and Boolean Algebras}
An important way of constructing a desired ultrafilter on $\kappa$ is
to use the construction of an ultrafilter $\E$ on the Boolean algebra
$\boB={\cal P}(\kappa)/\D$ for some filter $\D$ on $\kappa$ 
as an intermediate step. The construction of $\E$ has a great deal of 
flexibility when $\boB$ contains a large free (or $\kappa$-free)
subalgebra. In this paper we always construct an ultrafilter
$\E$ on $\boB$ first such that $\boul$, the Boolean ultrapower of $\omega$
modulo $\E$, has some desired properties, and then use $\E$ to define an 
ultrafilter $\F$ on $\kappa$ so that $\ul$, the ultrapower of $\omega$
modulo $\F$, is isomorphic to $\boul$. In each
case a large cardinal is used to construct $\D$ so that
$\boB={\cal P}(\kappa)/\D$ always contains a large free (or $\kappa$-free) 
subalgebra. 

\cite{Mansfield} is recommanded  for the general theory of 
the Boolean ultrapower of arbitrary models. Here we give 
the definitions and facts needed in this paper to keep it
self-contained. 

Throughout this paper we use $\kappa,\lambda,\eta$, etc. for
infinite cardinals, $\alpha,\beta,\gamma$, 
etc. for ordinals, and $k,l,m,n$, etc. for natural numbers. 
Let $A$ and $B$ be two sets. We denote by $A^B$ for
the set of all functions from $B$ to $A$ (except in the case when $B$ is 
a Boolean algebra). We also write $\kappa^\lambda$
for exponents in cardinal arithmetic, and this should be clear from
the context. Let $\kl$ be the set of all subsets of $\lambda$
of size $<\kappa$. Let $\D,\E,\F$, etc. denote filters or
ultrafilters, and let $\boB,\boC$, etc. denote Ba or cBa, {\em i.e.}
Boolean algebras or complete Boolean algebras.

We shall not distinguish a Ba $(\boB;\vee,\wedge,-,0,1)$
from its base set $\boB$. For any $S\subseteq\boB$ let $\bigvee S$
($\bigwedge S$) be the least upper bound (greatest lower bound) of
$S$ in $\boB$, provided it exists. By an anti-chain in $\boB$ we mean
a subset $A\subseteq\boB$ such that for any $a,b\in A$, $a\not=b$
implies $a\wedge b=0$. A Ba $\boB$ has $\kappa$-c.c. iff
every anti-chain $A$ in $\boB$ has size $<\kappa$. 
$\omega_1$-c.c. is called also c.c.c.

We write $\boC\subseteq\boB$ to denote that $\boC$ is a subalgebra of
$\boB$ and for any $S\subseteq \boC$, if $\bigvee S=c$ in $\boC$
and $\bigvee S=b$ in $\boB$, then $b=c$. Hence we shall not distinguish 
$\bigvee$ and $\bigwedge$ in $\boC$ or in $\boB$.
 
A $\boB$ is called a $\kappa$-complete Ba or a $\kappa$-cBa iff for any 
$S\subseteq\boB$, $|S|<\kappa$
implies $\bigvee S\in\boB$. $\boB$ is complete iff it is $\kappa$-complete
for every $\kappa$. $\omega_1$-complete is also called countably complete. 
Given $a\in\boB$, let $a^0$ denote $a$ and $a^1$ denote
$-a$. Let $\boC\subseteq\boB$. Then a sequence 
$\{a_{\alpha}:\alpha<\lambda\}$ in $\boB$ is called $\kappa$-independent
over $\boC$ iff for any $\sigma\in\kl$, for any $h\in 2^{\sigma}$ and
for any $c\in\boC\smallsetminus\{0\}$ one has
\[c\wedge (\bigwedge\{a^{h(\alpha)}_{\alpha}:\alpha\in\sigma\})\not=0.\]
A sequence in $\boB$ is $\kappa$-independent iff it is $\kappa$-independent
over $\{0,1\}$. A sequence $\{\boC_{\alpha}:\alpha<\lambda\}$ of 
subalgebras of $\boB$ is called $\kappa$-independent iff for any
$\sigma\in\kl$ and for any $a_{\alpha}\in\boC_{\alpha}\smallsetminus\{0\}$
\[\bigwedge_{\alpha\in\sigma}a_{\alpha}\not=0.\]
We shall omit $\kappa$ when $\kappa=\omega$.
Let $\boB$ be $\kappa$-complete. Then for any $S\subseteq\boB$ 
we denote by $\langle S\rangle_{\kappa}\subseteq\boB$ the 
$\kappa$-complete subalgebra
of $\boB$ generated by $S$. For any $\boB$ let $\bar{\boB}$ denote
the completion of $\boB$. A Ba $\boB$ is called 
$\kappa$-free iff there exists a $\kappa$-independent sequence 
$\{a_{\alpha}:\alpha<\lambda\}$ in $\boB$ such that
\[\boB=\langle\{a_{\alpha}:\alpha<\lambda\}\rangle_{\kappa}.\]
The cardinal $\lambda$ above is called the dimension of $\boB$.
Note that if $\kappa^{<\kappa}=\kappa$, then a $\kappa$-free Ba
has $\kappa^+$-c.c..
A Ba is free if it is $\omega$-free.
Given two Ba's $\boB\subseteq\boC$, a homomorphism\footnote{$r:\boC\mapsto
\boB$ is a homomorphism iff 
$r(a\vee b)=r(a)\vee r(b)$, and $r(-a)=-r(a)$. Hence one has
$r(a\wedge b)=r(a)\wedge r(b)$, $r(0)=0$ and $r(1)=1$. Note that
$r(\bigvee S)$ may not be same as $\bigvee\{r(a):a\in S\}$.}
$r:\boC\mapsto\boB$ is called a retraction iff $r\!\res\!\boB$ 
is an identity map.

Let $\boB$ be $\omega_1$-complete. Let
\[\bom=\{t:t\mbox{ is a function from }\omega\mbox{ to }\boB
\mbox{ such that }\]\[(1)\;\forall n\not=m\;(t(n)\wedge t(m)=0)\mbox{ and } 
(2)\;\bigvee\{t(n):n\in\omega\}=1.\}.\]
For any ultrafilter $\E$ on $\boB$ let $\sim_{\E}$ denote the equivalence
relation on $\bom$ such that $s\sim_{\E} t$ 
iff $\bigvee_{n\in\omega}(s(n)\wedge t(n))\in\E$ for any $s,t\in\bom$. 
For any $t\in\bom$
let $t_{\E}$ denote the equivalence class containing $t$.
Then $\boul$ is the set $\{t_{\E}:t\in\bom\}$ together with the total 
order $<_{\E}$, which is defined by letting $s_{\E}<_{\E}t_{\E}$ iff
$\bigvee_{m<n}(s(m)\wedge t(n))\in\E$. 
We shall write $<$ instead of $<_{\E}$ when the meaning is clear.
Note that if one identifies each $n\in\omega$ with $(t_n)_{\E}$,
where $t_n(n)=1$ and $t_n(m)=0$ for any $m\not=n$, then
$\boul$ is just an end-extension of $\omega$. To make it more intuitive
we often write $[\![ s<t ]\!]$ and $[\![ s=t ]\!]$
for the terms $\bigvee_{m<n}(s(m)\wedge t(n))$ and 
$\bigvee_{n\in\omega}(s(n)\wedge t(n))$, respectively.

Suppose $\D$ is a filter on $\kappa$ and let $I$ be the dual ideal of $\D$.
We often write ${\cal P}(\kappa)/{\D}$ instead of ${\cal P}(\kappa)/I$,
the quotient Boolean algebra of ${\cal P}(\kappa)$ modulo $I$.
For any $A\subseteq\kappa$ let $[A]_{\D}$ denote the equivalence
class of $A$ in ${\cal P}(\kappa)/\D$.
 
Suppose $\D$ is a countably complete filter on $\kappa$. Then the Boolean 
algebra $\boB={\cal P}(\kappa)/\D$ is countably complete, and for any
$A_n\subseteq\kappa$ one has 
\[[\bigcup_{n\in\omega}A_n]_{\D}=\bigvee_{n\in\omega}[A_n]_{\D}.\]
Let $\E$ be an ultrafilter on $\boB$ and let
\[\F=\{A\subseteq\kappa:[A]_{\D}\in\E\}.\]
Then it is easy to see that $\F$ is an ultrafilter on $\kappa$ extending
$\D$. We want to show that $\boul$ and $\ul$ are isomorphic.

For any $f\in\omega^{\kappa}$ let $\hat{f}$ be a function from $\omega$
to $\boB$ such that $\hat{f}(n)=[f^{-1}(n)]_{\D}$. Then it is easy to 
check that $\hat{f}\in\bom$. Let $i$ be the map from $\ul$ to $\boul$
such that $i(f_{\F})=\hat{f}_{\E}$.

\begin{Lemma} \label{iso} (folklore)

The map $i$ is an isomorphism from $\ul$ to $\boul$.

\end{Lemma}

\noindent {\bf Proof:}\quad
We show first that $i$ is a well-defined one-one function,
which preserves the order.
This is because that for any two functions $f,g\in\omega^{\kappa}$, one has  
\[f_{\F}<g_{\F}\Leftrightarrow
\{\alpha<\kappa:f(\alpha)<g(\alpha)\}\in\F\Leftrightarrow\] 
\[\bigcup_{m<n}(f^{-1}(m)\cap g^{-1}(n))\in\F\Leftrightarrow
[\bigcup_{m<n}(f^{-1}(m)\cap g^{-1}(n))]_{\D}\in\E.\]
But one has also
\[[\bigcup_{m<n}(f^{-1}(m)\cap g^{-1}(n))]_{\D}=
\bigvee_{m<n}([f^{-1}(m)]_{\D}\wedge [g^{-1}(n)]_{\D})=
\bigvee_{m<n}(\hat{f}(m)\wedge \hat{g}(n)).\]
So it is true that
\[f_{\F}<g_{\F}\Leftrightarrow\hat{f}_{\E}<\hat{g}_{\E}\Leftrightarrow 
i(f_{\F})<i(g_{\F}).\]

Then we show that $i$ is onto. Given any $t\in\bom$, by countable 
completeness we could choose $A_n\subseteq\kappa$ inductively such that
$\{A_n:n\in\omega\}$ is a partition of $\kappa$ and $[A_n]_{\D}=t(n)$.
Define an $f\in\omega^{\kappa}$ such that $f(\alpha)=n$ iff 
$\alpha\in A_n$. Clearly, $i(f_{\F})=t_{\E}$. \quad $\Box$

\begin{Remark}

By this lemma to find an ultrafilter
$\F$ on $\kappa$ such that $\ul$ has some desired properties, 
it suffices to find a countably complete filter $\D$ on 
$\kappa$ and an ultrafilter
$\E$ on $\boB={\cal P}(\kappa)/\D$ such that $\boul$ has such properties.

\end{Remark}
 
\section{Strong Limit Cardinal as the Size of an Ultrapower of $\omega$}
An ultrafilter $\F$ on $\kappa$ is called uniform iff $|A|=\kappa$ for every
$A\in\F$. An ultrafilter $\F$ on $\kappa$ is called regular iff there exists
a family $\C\subseteq\F$ such that $|\C|=\kappa$ and for any infinite 
subfamily $\C_0\subseteq\C$, $\bigcap\C_0=\emptyset$.

The study of the cardinality of ultrapowers started in the 1960's.
The subject is very closely related to the regularity of ultrafilters. 
Regular ultrafilters were introduced in \cite{FMSc} and \cite{Keisler1}. 
It was proved there that if $\F$ is a regular ultrafilter 
on $\kappa$, then $\ul$ has size
$2^{\kappa}$. So it is natural to ask whether it is possible to
have an ultrafilter $\F$ on $\kappa$ such that the size of
$\ul$ does not have the form $2^{\kappa}$ for any $\kappa$. In fact,
Chang and Keisler asked this in the following forms 
(see page 252 of \cite{CK}).

\begin{Question} \label{question1}

Is it possible that the cardinality of an ultrapower of $\omega$ 
is a singular (strong limit) cardinal?

\end{Question}

\begin{Question}\label{question2} 

Is it possible that the cardinality of an ultrapower of $\omega$ 
is a strongly inaccessible cardinal? 

\end{Question}

The original form of Question \ref{question1} in \cite{CK} has no requirement 
for the cardinality to be a strong limit cardinal.
Since the cardinal $2^{\kappa}$ could be singular, we would like to 
make it more specific by requiring the singular cardinality be also 
a strong limit.

Obviously, a positive answer to either Question \ref{question1} or 
Question \ref{question2} would imply the existence of a uniform
non-regular ultrafilter. Since the existence of a uniform non-regular
ultrafilter was unclear in the early 1970's, when the first edition of
\cite{CK} was published, people were more interested
in a general question of Chang and Keisler concerning the existence
of uniform non-regular ultrafilters. A lot of work has been done
since then for solving the general question. On one hand,
Prikry \cite{Prikry}, Ketonen \cite{Ketonen}, Donder \cite{Donder},
etc. showed that one may need to assume the consistency of some large
cardinals to construct a uniform non-regular ultrafilter. For example,
Donder proved that if there is no inner model containing a measurable
cardinal, then (1) every uniform ultrafilter on a singular cardinal
is regular, (2) every uniform ultrafilter on a regular cardinal $\kappa$
with $(\kappa^+)^K=\kappa^+$ is regular, where $K$ is the 
Dodd-Jensen core model.
On the other hand, Magidor \cite{Magidor}, Laver \cite{Laver2}, 
Foreman, Magidor and Shelah \cite{FMSh}, etc. showed that it is possible to
construct a uniform non-regular ultrafilter with the help of some large
cardinal axioms. For example, in \cite{FMSh} it is proved that,
assuming $\kappa$ is a huge cardinal and $\mu<\kappa$ is a regular cardinal,
there is a forcing extension preserving every cardinal $\leqslant\mu$,
in which there exists a uniform (fully) non-regular ultrafilter on $\mu^+$.
However, as far as we know, neither Question \ref{question1} nor 
Question \ref{question2} has been answered yet. In this section
we are going to give positive answers to both questions by assuming
the consistency of a supercompact cardinal. First, we need to introduce
the Laver-indestructibility of a supercompact cardinal. See \cite{Jech}
or \cite{KRS} for the basic facts of supercompact cardinals

A supercompact cardinal $\kappa$ is called Laver-indestructible iff
$\kappa$ remains supercompact in any $\kappa$-directed closed forcing
extension. The reader should consult \cite{Laver1} to see how one
makes a supercompact cardinal Laver-indestructible by a $\kappa$-c.c. 
forcing of size $\kappa$.

\begin{Theorem} \label{singular}

Suppose $\kappa$ is a Laver-indestructible supercompact cardinal.
Let $\eta\geqslant\kappa$ be such that $\eta^{<\kappa}=\eta$ and let
$\lambda\leqslant 2^{\eta}$ be such that $\lambda^{\kappa}=\lambda$. 
Then there exists an ultrafilter $\F$ 
on $\eta$ such that $|\omega^{\eta}/\F|=\lambda$.

\end{Theorem}

\begin{Remark}

Note that $\lambda$ could be a singular strong limit cardinal, say,
$\lambda=\eta=\beth_{\kappa^+}(\kappa)$. So Theorem \ref{singular} answers 
Question \ref{question1}. If one assumes that
there is a strongly inaccessible cardinal $\lambda$ above $\kappa$, then 
Question \ref{question2} could also be answered. But the next theorem 
shows that this extra assumption is not necessary.

\end{Remark}

\begin{Theorem} \label{inaccessible}

Suppose $\kappa$ is a Laver-indestructible supercompact cardinal.
Then there exists an ultrafilter $\F$ on $\kappa$ such that
$|\ul|=\kappa$.

\end{Theorem}

Let $\kappa$ be a strongly compact cardinal and $\lambda\geqslant\kappa$.
Then it is well-known (Lemma 33.1 of \cite{Jech}) that any 
$\kappa$-complete filter on $\lambda$ could be extended to a
$\kappa$-complete ultrafilter on $\lambda$. Note that a supercompact
cardinal is strongly compact.

\begin{Lemma} \label{key1}

Suppose $\kappa$ is a Laver-indestructible supercompact cardinal.
Let $\eta\geqslant\kappa$ be such that $\eta^{<\kappa}=\eta$ and let
$\kappa\leqslant\lambda\leqslant 2^{\eta}$.
Then there exists a $\kappa$-complete filter $\D$ on $\eta$ such that
$\boB={\cal P}(\eta)/\D$ has a $\kappa$-free dense subalgebra 
$\boC\subseteq\boB$ of dimension $\lambda$.

\end{Lemma}

\noindent {\bf Proof:}\quad Let ${\Bbb P}=Fn(\lambda,2,\kappa)$ 
(see page 211 of \cite{Kunen} for the definition) be the forcing
notion for adding $\lambda$ Cohen subsets of $\kappa$. Note that ${\Bbb P}$
is $\kappa$-directed closed. By Laver-indestructibility $\kappa$ is still
supercompact in $V^{\Bbb P}$. Suppose $G\subseteq {\Bbb P}$ is a
$V$-generic filter and $g=\bigcup G$. Then $g$ is a function from
$\lambda$ to $2$ in $V[G]$. In $V$ one can choose a sequence
$\{A_{\alpha}:\alpha<\lambda\}$ of subsets of $\eta$ such that for any 
$\sigma\in\kl$ and for any $h\in 2^{\sigma}$ one has
\[|\bigcap_{\alpha\in\sigma}A^{h(\alpha)}_{\alpha}|=\eta,\]
where we denote $A^0=A$ and $A^1=\eta\smallsetminus A$ for 
$A\subseteq\eta$. The existence of such a sequence is guaranteed by
$\eta^{<\kappa}=\eta$ (see page 288 of \cite{Kunen}). 
In $V^{\Bbb P}$ the set 
$\{A^{g(\alpha)}_{\alpha}:\alpha<\lambda\}$ forms a 
$\kappa$-complete filter base. Hence there exists an $\kappa$-complete
ultrafilter ${\cal H}_g$ on $\eta$ extending 
$\{A^{g(\alpha)}_{\alpha}:\alpha<\lambda\}$. Now back in $V$ we define
\[\D=\{A\subseteq\eta:\;\Vdash A\in {\cal H}_{\dot{g}}\}.\]
It is clear that $\D$ is a $\kappa$-complete filter in $V$.
Let $I$ be the dual ideal of $\D$.

\medskip

{\bf Claim \ref{key1}.1.}\quad For any $\sigma\in\kl$ 
and for any $h\in 2^{\sigma}$
one has $\bigcap_{\alpha\in\sigma}A^{h(\alpha)}_{\alpha}\not\in I$.

\medskip

Proof of Claim \ref{key1}.1:\quad Since for any $\alpha\in\sigma$
\[h\Vdash A^{h(\alpha)}_{\alpha}\in {\cal H}_{\dot{g}},\]
Then one has 
\[h\Vdash\bigcap_{\alpha\in\sigma}A^{h(\alpha)}_{\alpha}
\in {\cal H}_{\dot{g}}.\] Hence one has
\[\not\Vdash\eta\smallsetminus
\bigcap_{\alpha\in\sigma}A^{h(\alpha)}_{\alpha}
\in {\cal H}_{\dot{g}}.\quad \Box (\mbox{Claim \ref{key1}.1})\]

{\bf Claim \ref{key1}.2.}\quad If $A\not\in I$, then there exists a $\sigma\in
\kl$ and an $h\in 2^{\sigma}$ such that $\bigcap_{\alpha\in\sigma}
A^{h(\alpha)}_{\alpha}\smallsetminus A\in I$.

\medskip

Proof of Claim \ref{key1}.2:\quad Suppose $A\not\in I$. Then 
\[\not\Vdash\eta\smallsetminus A\in {\cal H}_{\dot{g}}.\]
So there exists an $h\in {\Bbb P}$ such that
\[h\Vdash\eta\smallsetminus A\not\in {\cal H}_{\dot{g}}.\]
This means that \[h\Vdash A\in {\cal H}_{\dot{g}}.\]
We now want to show that $\bigcap_{\alpha\in\sigma}
A^{h(\alpha)}_{\alpha}\smallsetminus A\in I$, where $\sigma=dom(h)$.
Suppose not. Then 
\[\not\Vdash\eta\smallsetminus (\bigcap_{\alpha\in\sigma}
A^{h(\alpha)}_{\alpha}\smallsetminus A)\in {\cal H}_{\dot{g}}.\]
Hence there exists $h'\in {\Bbb P}$ such that
\[h'\Vdash\bigcap_{\alpha\in\sigma} A^{h(\alpha)}_{\alpha}\smallsetminus A
\in {\cal H}_{\dot{g}}.\] This implies
\[h'\Vdash\bigcap_{\alpha\in\sigma} A^{h(\alpha)}_{\alpha}\in 
{\cal H}_{\dot{g}}\mbox{ and }h'\Vdash A\not\in {\cal H}_{\dot{g}}.\] 
Let $\sigma'=dom(h')$. For any $\alpha\in\sigma'\cap\sigma$ one has
\[h'\Vdash A^{h(\alpha)}_{\alpha}\in {\cal H}_{\dot{g}}\mbox{ and }
h'\Vdash A^{h'(\alpha)}_{\alpha}\in {\cal H}_{\dot{g}}.\]
Hence $h(\alpha)=h'(\alpha)$. So $h\cup h'\in {\Bbb P}$.
But this contradicts \[h\Vdash A\in {\cal H}_{\dot{g}}.
\quad \Box (\mbox{Claim \ref{key1}.2})\]

Let $a_{\alpha}=[A_{\alpha}]_{\D}$. By Claim \ref{key1}.1 the sequence
$\{a_{\alpha}:\alpha<\lambda\}$ is $\kappa$-independent. Let 
$\boC=\langle\{a_{\alpha}:\alpha<\lambda\}\rangle_{\kappa}$.
Then by Claim \ref{key1}.2 $\boC$ is dense in $\boB={\cal P}(\eta)/\D$.
\quad $\Box$

\bigskip

Now we are ready to prove the theorems.

\bigskip

\noindent {\bf Proof of Theorem \ref{singular}:}\quad
Let $\lambda\leqslant 2^{\eta}$ be such that 
$\lambda^{\kappa}=\lambda$. Let $\D$
be a $\kappa$-complete filter on $\eta$ obtained in Lemma \ref{key1}
such that the Ba $\boB= {\cal P}(\eta)/\D$ has a
$\kappa$-free dense subalgebra $\boC\subseteq\boB$ of dimension $\lambda$.
It is clear that $|\boC|=\lambda$ and every element $a$ in $\boB$ is the
least upper bound of an anti-chain in $\boC$. Since $\boC$ has 
$\kappa^+$-c.c., every anti-chain in $\boC$ has size $\leqslant\kappa$. 
Hence $|\boB|\leqslant |\boC|^{\kappa}=\lambda^{\kappa}=
\lambda$. This shows that for any ultrafilter $\E$ on $\boB$ one has
$|\boul|\leqslant\lambda^{\omega}=\lambda$.

We now want to show that there exists an ultrafilter $\E$ on $\boB$ such that
$|\boul|\geqslant\lambda$. Let $\{a_{\alpha,n}:\alpha<\lambda,\;n\in\omega\}$ 
be a $\kappa$-independent sequence in $\boC$ such that 
\[\boC=\langle\{a_{\alpha,n}:\alpha<\lambda,n\in\omega\}\rangle_{\kappa}.\]
For each $\alpha<\lambda$ let $t_{\alpha}\in\bom$ be such that
\[t_{\alpha}(0)=-(\bigvee_{n\in\omega}a_{\alpha,n})\;\mbox{ and }\;
t_{\alpha}(n+1)=a_{\alpha,n}\wedge -(\bigvee_{m<n}a_{\alpha,m}).\]
Let $\E_0=\{[\![t_{\alpha}<t_{\beta}]\!]:\alpha<\beta<\lambda\}$.

\medskip

{\bf Claim \ref{singular}.1.}\quad $\E_0$ has the finite intersection property.

\medskip

Proof of Claim \ref{singular}.1:\quad
Since $[\![t_{\alpha}< t_{\beta}]\!]\wedge 
[\![t_{\beta}< t_{\gamma}]\!]
\leqslant [\![t_{\alpha}< t_{\gamma}]\!]$, then one needs only to
show that for any $\alpha_0<\alpha_1<\ldots<\alpha_k$ in $\lambda$
\[\bigwedge_{n<k}[\![t_{\alpha_n}<t_{\alpha_{n+1}}]\!]\not=0.\]
Choose any $m_0<m_1<\ldots<m_k$ in $\omega$. Then one has
\[\bigwedge_{n\leqslant k}t_{\alpha_n}(m_n+1)=
\bigwedge_{n<k}(t_{\alpha_n}(m_n+1)\wedge t_{\alpha_{n+1}}(m_{n+1}+1))
\leqslant \bigwedge_{n<k}[\![t_{\alpha_n}<t_{\alpha_{n+1}}]\!].\]
But one has also that
\[\bigwedge_{n\leqslant k}t_{\alpha_n}(m_n+1)=
\bigwedge_{n\leqslant k}(a_{\alpha_n,m_n}\wedge (\bigwedge_{l<m_n}
(-a_{\alpha_n,l})))\not=0\] by the independence of $a_{\alpha,n}$'s.
\quad $\Box$ (Claim \ref{singular}.1)

\medskip

Let $\E$ be an ultrafilter on $\boB$ extending $\E_0$. Then there is a 
strictly increasing sequence $\{(t_{\alpha})_{\E}:\alpha<\lambda\}$
in $\boul$. Hence $|\boul|\geqslant\lambda$.
Now the theorem follows from Lemma \ref{iso}.\quad $\Box$

\bigskip

\noindent {\bf Proof of Theorem \ref{inaccessible}:}\quad
Again by Lemma \ref{key1}\footnote{The full strength of
supercompactness is not necessary here; for example, 
$2^{\kappa}$-supercompactness suffices.} 
one can find a $\kappa$-complete filter $\D$
such that $\boB= {\cal P}(\kappa)/\D$ has a $\kappa$-free, $\kappa^+$-c.c.
dense subalgebra $\boC\subseteq\boB$ of dimension $\kappa$ generated by a
$\kappa$-independent sequence $\{a_{\alpha,n}:\alpha<\kappa,n\in\omega\}$.
Let $t_{\alpha}\in\bom$ be same as in the proof of Theorem \ref{singular}
for every $\alpha<\kappa$.
For any successor ordinal $\alpha$ let 
$\boC_{\alpha}=\langle\{a_{\beta,n}:\beta<\alpha,n\in\omega\}\rangle_{\kappa}$
and for any limit ordinal $\alpha<\kappa$ let 
$\boC_{\alpha}=\bigcup_{\beta<\alpha}\boC_{\beta}$.
Note that for any successor $\alpha$ $\boC_{\alpha}$ is
atomic, for any limit $\alpha$ $\boC_{\alpha}$ is not 
$\kappa$-complete, and for any $\alpha<\kappa$ 
$|\boC_{\alpha}|<\kappa$. It is easy to see that 
$\boC=\bigcup_{\alpha<\kappa}\boC_{\alpha}$.
For each $\alpha<\kappa$ let 
\[D_{\alpha}=\{\bigvee A:A\mbox{ is a maximal anti-chain in }\boC_{\alpha}\}.\] 

\medskip

{\bf Claim \ref{inaccessible}.1.}\quad The set $(\bigcup_{\alpha<\kappa}
D_{\alpha})\cup\{[\![t_{\alpha}<t_{\beta}]\!]:\alpha<\beta<\kappa\}$
has the finite intersection property.

\medskip

Proof of Claim \ref{inaccessible}.1:\quad
We show by induction on $\gamma$ that the set
\[(\bigcup_{\alpha<\gamma}D_{\alpha})\cup\{[\![t_{\alpha}<t_{\beta}]\!]:
\alpha<\beta<\gamma\}\] has the finite intersection property.
This is trivial when $\gamma=0$ or $\gamma$ is a limit ordinal. 
Let's assume that $\gamma=\beta+1$ for some $\beta<\kappa$.
It suffices to show that for any $a\in\bigcup_{\alpha<\beta}\boC_{\alpha}
\smallsetminus\{0\}$, for any set of maximal anti-chains 
$A_0,A_1,\ldots,A_k\subseteq\boC_{\beta}$ and for any $\alpha<\beta$
\[a\wedge (\bigwedge_{n\leqslant k}\bigvee A_n)\wedge [\![t_{\alpha}<
t_{\beta}]\!]\not=0.\]
Since $\bigvee_{n\in\omega}t_{\alpha}(n)=1$ there is an $m\in\omega$
such that $a\wedge t_{\alpha}(m)\not=0$. It is easy to see that
$a\wedge t_{\alpha}(m)\wedge t_{\beta}(m+1)\not=0$ because
$t_{\beta}(m+1)$ is independent over $\boC_{\beta}$ and
$a\wedge t_{\alpha}(m)\in\boC_{\beta}$.
Now using the maximality of $A_n$'s one could find $a_n\in A_n$ for each
$n\leqslant k$ inductively on $n$ such that 
\[a\wedge t_{\alpha}(m)\wedge t_{\beta}(m+1)\wedge
\bigwedge_{n\leqslant k}a_n\not=0.\]
This finishes the proof because
\[a\wedge t_{\alpha}(m)\wedge t_{\beta}(m+1)\wedge\bigwedge_{n\leqslant k}a_n
\leqslant a\wedge (\bigwedge_{n\leqslant k}\bigvee A_n)\wedge 
[\![t_{\alpha}< t_{\beta}]\!].\]
\quad $\Box$(Claim \ref{inaccessible}.1)

\medskip

By Claim \ref{inaccessible}.1 one can find an ultrafilter $\E$ on 
$\boB$ extending the set 
\[(\bigcup_{\alpha<\kappa}
D_{\alpha})\cup\{[\![t_{\alpha}<t_{\beta}]\!]:\alpha<\beta<\kappa\}.\]
Clearly, $\{(t_{\alpha})_{\E}:\alpha<\kappa\}$ is a strictly increasing
sequence of length $\kappa$ in $\boul$. Hence $|\boul|\geqslant\kappa$.
We need to show that $|\boul|\leqslant\kappa$.

\medskip

{\bf Claim \ref{inaccessible}.2.}\quad
For any maximal anti-chain $\{b_{\alpha}:\alpha<\kappa\}$ in $\boC$ there 
exists a $\delta<\kappa$ such that $\bigvee\{b_{\alpha}:\alpha<\delta\}\in\E$.

\medskip

Proof of Claim \ref{inaccessible}.2:\quad
Note that $|\boC_{\alpha}|<\kappa$ for any $\alpha<\kappa$. 
Using the inaccessibility of $\kappa$ one could show that 
there exists a $\delta<\kappa$
such that $\{b_{\alpha}:\alpha<\delta\}$ is a maximal anti-chain in
$\boC_{\delta}$. Hence $\bigvee\{b_{\alpha}:\alpha<\delta\}\in\E$.
\quad $\Box$(Claim \ref{inaccessible}.2)

\medskip

{\bf Claim \ref{inaccessible}.3.}\quad For any $t\in\bom$ there exists
an $s\in\omega^{\boC}$ such that $t_{\E}=s_{\E}$.

\medskip

Proof of Claim \ref{inaccessible}.3:\quad
Since $\boC$ is dense in $\boB$ and $\boC$ has $\kappa^+$-c.c., 
then there exists a maximal anti-chain $\{b_{\alpha}:\alpha<\kappa\}$
in $\boC$, which refines $\{t(n):n\in\omega\}$, {\em i.e.} for any
$\alpha<\kappa$ there exists an $n\in\omega$ such that $b_{\alpha}\leqslant
t(n)$. By Claim \ref{inaccessible}.2 there is a $\delta<\kappa$
such that $\bigvee\{b_{\alpha}:\alpha<\delta\}\in\E$. Define 
$s\in\omega^{\boC}$ such that 
\[s(0)=(\bigvee\{b_{\alpha}:\alpha<\delta,\;b_{\alpha}\leqslant t(0)\})
\vee (-\bigvee\{b_{\alpha}:\alpha<\delta\})\] and
\[s(n+1)=\bigvee\{b_{\alpha}:\alpha<\delta,\;b_{\alpha}\leqslant t(n+1)\}\]
for every $n\in\omega$.
Here we use the fact that $\boC$ is $\kappa$-complete and $\delta<\kappa$.
It is now easy to check that
\[[\![t=s]\!]=\bigvee_{n\in\omega}(t(n)\wedge s(n))\geqslant
\bigvee\{b_{\alpha}:\alpha<\delta\}\in\E.\]
Hence $|\boul|\leqslant |\boC|^{\omega}=\kappa$.
Now the theorem follows from Lemma \ref{iso}.
\quad $\Box$

\section{$\lambda$-Archimedean Ultrapowers}
Let $\cal L$ be any first-order language including 
symbols for number theory, say, ${\Bbb N},+,\cdot,<$ and $S$. An 
$\cal L$-structure $\frak A$ is called a (standard or nonstandard) 
model of PA (PA stands for Peano Arithmetic) iff 
$({\cal N}^{\frak A};+^{\frak A},\cdot^{\frak A},<^{\frak A},S^{\frak A})$ 
is a (standard or nonstandard) model of PA (respectively).
Given a cardinal $\lambda$, a model $\frak A$ of PA is called 
$\lambda$-Archimedean iff $|{\Bbb N}^{\frak A}|=\lambda$ and
for every $n\in{\Bbb N}^{\frak A}$, $|\{0,1,\ldots,n\}^{\frak A}|<\lambda$.

In \cite{KS} Keisler and Schmerl asked whether one could find an ultrapower
of a standard model of PA, which is $\lambda$-Archimedean for some cardinal
$\lambda>\omega$. It is clear that the question remains same if one replaces the
standard model of PA by $(\omega;<)$. If $\F$ is a regular ultrafilter on
$\kappa$, then $\ul$ will never be 
$\lambda$-Archimedean for any $\lambda>\omega$ because there exists
an $x\in\ul$ such that $|\{y\in\ul:y<x\}|=2^{\kappa}=|\ul|$ (see, for
example, \cite{Koppelberg}).
So a positive answer to the question implies the existence of a uniform
non-regular ultrafilter.

In this section we construct several $\lambda$-Archimedean ultrapowers
under some large cardinal assumptions.

First, we list some lemmas needed in the proof of theorems.  

\begin{Lemma}\label{Sikorski} (Sikorski) 

Let $\boB$ be a cBa and $\boB\subseteq\boC$. Then there is 
a retraction $r$ from $\boC$ to $\boB$. Furthermore, if
$I$ is an ideal of $\boC$ and $I\cap\boB=\{0\}$, then
one could require that $r(a)=0$ for every $a\in I$.

\end{Lemma}

\noindent {\bf Proof:}\quad See page 70 of \cite{MB} for the first 
assertion. For the second assertion one considers the fact that
$\boB$ could be viewed as a subalgebra of $\boC/I$. Then one uses
the retraction $\bar{r}$ from $\boC/I$ to $\boB$ to induce a
retraction $r$ from $\boC$ to $\boB$.\quad $\Box$

\begin{Lemma} \label{Mansfield} (folklore)

Let $\boB$ and $\boC$ be two $\omega_1$-cBa's such that $\boB\subseteq\boC$. 
Let $\E$ be an ultrafilter on $\boB$ and $\E'$ be
an ultrafilter on $\boC$ such that $\E\subseteq\E'$. Then the inclusion
map $i$ from $\boul$ to $\omega^{\boC}/\E'$ such that
$i(t_{\E})=t_{\E'}$ is an elementary embedding.

\end{Lemma}

See page 576 of \cite{Koppelberg} for a proof.
\relax From now on we shall view $\boul$ as a subset of $\omega^{\boC}/\E'$
via the embedding $i$ whenever $\E\subseteq\E'$.

\begin{Lemma} \label{Koppelberg} (Exercise IV.3.35 of \cite{Shelah2} or Lemma 1 of \cite{Koppelberg})

Let $\boB$ and $\boC$ be two $\omega_1$-cBa's such that $\boB\subseteq\boC$. 
Suppose $\E$ is an ultrafilter on $\boB$, 
$r:\boC\mapsto\boB$ is a retraction and $\E'=r^{-1}(\E)$. Then
$\boul$ is an initial segment of $\omega^{\boC}/\E'$. Furthermore,
for any $t\in\omega^{\boC}$ $t_{\E'}\in\boul$ iff $\bigvee_{n\in\omega}
r(t(n))\in\E$.

\end{Lemma}

%
%
%
%
%

\begin{Lemma}\label{Shelah1}

Let $\{\boB_{\alpha}:\alpha<\delta\}$ be a sequence of cBa's for some
limit ordinal $\delta$ and $\boB_{\delta},\boC$ be two Ba's. 
Suppose $\boB_{\alpha}\subseteq\boB_{\beta}
\subseteq\boB_{\delta}\subseteq\boC$ and $r_{\alpha,\gamma}:
\boB_{\gamma}\mapsto\boB_{\alpha}$ are
retractions such that $r_{\alpha,\beta}\circ r_{\beta,\gamma}=
r_{\alpha,\gamma}$ for any $\alpha<\beta<\gamma\leqslant\delta$. 
Given $x\in\boC$, let $y_{\alpha}\in\boB_{\alpha}$ be such that
\[\bigvee\{r_{\alpha,\delta}(a):a\leqslant x,\;a\in\boB_{\delta}\}\leqslant
y_{\alpha}\leqslant\bigwedge\{r_{\alpha,\delta}(a):a\geqslant x,
\;a\in\boB_{\delta}\}\] and $y_{\alpha}=r_{\alpha,\beta}(y_{\beta})$ for
any $\alpha<\beta<\delta$. Then there exist
retractions $p_{\alpha}:\langle\boB_{\delta}\cup\{x\}\rangle\mapsto
\boB_{\alpha}$ such that $p_{\alpha}\!\res\!\boB_{\delta}=
r_{\alpha,\delta}$, $p_{\alpha}=r_{\alpha,\beta}\circ p_{\beta}$
and $p_{\alpha}(x)=y_{\alpha}$ for any $\alpha<\beta<\delta$.

\end{Lemma}

\noindent {\bf Proof:}\quad Trivial.

\begin{Lemma}\label{Shelah2}

Let $\{\boB_{\alpha}:\alpha<\delta\}$, $\boB_{\delta}$, $\boC$,
$r_{\alpha,\gamma}$ for any $\alpha<\gamma\leqslant\delta$ be as in 
Lemma \ref{Shelah1}. Then there exist retractions 
$p_{\alpha}:\boC\mapsto\boB_{\alpha}$ such that
$p_{\alpha}\!\res\!\boB_{\delta}=r_{\alpha,\delta}$ and 
$p_{\alpha}=r_{\alpha,\beta}\circ p_{\beta}$ for any $\alpha<\beta<\delta$.

\end{Lemma}

\noindent {\bf Proof:}\quad
By Lemma \ref{Shelah1} and Zorn's Lemma, it suffices to show that
for any $x\in\boC$ there exists $y_{\alpha}\in\boB_{\alpha}$ such that
\[\bigvee\{r_{\alpha,\delta}(a):a\leqslant x,\;a\in\boB_{\delta}\}\leqslant
y_{\alpha}\leqslant\bigwedge\{r_{\alpha,\delta}(a):a\geqslant x,
\;a\in\boB_{\delta}\}\] 
and $y_{\alpha}=r_{\alpha,\beta}(y_{\beta})$ for
any $\alpha<\beta<\delta$. Given any $\alpha<\delta$, let
\[u_{\alpha}=\bigvee\{r_{\alpha,\delta}(a):
a\leqslant x,\;a\in\boB_{\delta}\}\]
and \[v_{\alpha}=\bigwedge\{r_{\alpha,\delta}(a):a\geqslant x,
\;a\in\boB_{\delta}\}\]
for each $\alpha<\delta$. It is easy to check
that $u_{\alpha}\leqslant r_{\alpha,\beta}(u_{\beta})\leqslant
r_{\alpha,\beta}(v_{\beta})\leqslant v_{\alpha}$ 
for any $\alpha<\beta<\delta$. Let 
$\cal S$ be the set of all sequences $\{\langle a_{\alpha},b_{\alpha}\rangle:
\alpha<\delta\}$ with $u_{\alpha}\leqslant a_{\alpha}\leqslant b_{\alpha}
\leqslant v_{\alpha}$ and $a_{\alpha}\leqslant r_{\alpha,\beta}(a_{\beta})
\leqslant r_{\alpha,\beta}(b_{\beta})\leqslant b_{\alpha}$ for any
$\alpha<\beta<\delta$. Define a partial order $\leqslant_{\cal S}$ on
$\cal S$ such that 
\[\{\langle a_{\alpha},b_{\alpha}\rangle:\alpha<\delta\}
\leqslant_{\cal S}\{\langle a'_{\alpha},b'_{\alpha}\rangle:\alpha<\delta\}\]
iff $a_{\alpha}\leqslant a'_{\alpha}$ and $b'_{\alpha}\leqslant b_{\alpha}$
for every $\alpha<\delta$. Suppose 
\[{\cal T}=\{\{\langle a^i_{\alpha},b^i_{\alpha}\rangle:
\alpha<\delta\}:i\in I\}.\] is a totally ordered subset
of $\cal S$. Let $a^{\cal T}_{\alpha}=\bigvee\{a^i_{\alpha}:i\in I\}$
and $b^{\cal T}_{\alpha}=\bigwedge\{b^i_{\alpha}:i\in I\}$
for each $\alpha<\delta$. It is easy to check that the sequence
$\{\langle a^{\cal T}_{\alpha},b^{\cal T}_{\alpha}\rangle:\alpha<\delta\}$
is in $\cal S$ and is an upper bound of $\cal T$. By Zorn's Lemma
there is a maximal element $\{\langle c_{\alpha},d_{\alpha}\rangle:
\alpha<\delta\}$ in $\cal S$.

\medskip

{\bf Claim \ref{Shelah2}.1.}\quad $r_{\alpha,\beta}(c_{\beta})=c_{\alpha}$
and $r_{\alpha,\beta}(d_{\beta})=d_{\alpha}$ for any $\alpha<\beta<\delta$.

\medskip

Proof of Claim \ref{Shelah2}.1:\quad
Suppose not. Without loss of generality let $\beta_0<\delta$ be the 
smallest such that there exists an $\alpha_0<\beta_0$ with
$c_{\alpha_0}<r_{\alpha_0,\beta_0}(c_{\beta_0})$.
Now let $\bar{c}_{\alpha}=c_{\alpha}$ when $\alpha\geqslant\beta_0$
and $\bar{c}_{\alpha}=r_{\alpha,\beta_0}(c_{\beta_0})$ when $\alpha<\beta_0$.
Then it is easy to check that the sequence
$\{\langle\bar{c}_{\alpha},d_{\alpha}\rangle:\alpha<\delta\}$ is in
$\cal S$ and is strictly greater than
$\{\langle c_{\alpha},d_{\alpha}\rangle:\alpha<\delta\}$.
This contradicts the maximality of $\{\langle c_{\alpha},d_{\alpha}\rangle:
\alpha<\delta\}$.\quad $\Box$(Claim \ref{Shelah2}.1)

\medskip

Clearly the sequence $\{\langle c_{\alpha},d_{\alpha}\rangle:\alpha<\delta\}$
is what we want.\quad $\Box$

\begin{Theorem}\label{arch1} 

Assume $\kappa$ is a Laver-indestructible supercompact
cardinal and $\eta>\kappa$ such that $\eta^{<\kappa}=\eta$. 
Suppose $\lambda\leqslant 2^{\eta}$ such that
$\theta^{\kappa}<\lambda$ for any $\theta<\lambda$ and 
$\lambda^{\kappa}=\lambda$.
Then there exists an ultrafilter $\F$ on $\eta$ such that
$\omega^{\eta}/\F$ is $\lambda$-Archimedean.

\end{Theorem}

\noindent {\bf Proof:}\quad
By Lemma \ref{key1} there exists a $\kappa$-complete filter $\D$
on $\eta$ such that $\boB={\cal P}(\eta)/\D$ contains a $\kappa$-free
dense subalgebra $\boC$ of dimension $\lambda$. Without loss of generality
let \[\boC=\langle\{a_{\alpha,n}:\alpha<\lambda,\;n\in\omega\}
\rangle_{\kappa},\] where $\{a_{\alpha,n}:\alpha<\lambda,\;n\in\omega\}$
is a $\kappa$-independent sequence in $\boC$.
For each $\beta\leqslant\lambda$ let
\[\boC_{\beta}=\langle\{a_{\alpha,n}:\alpha<\beta,
\;n\in\omega\}\rangle_{\kappa}\]
and let $\boB_{\beta}=\bar{\boC}_{\beta}$.
Note that $\boC_0=\boB_0=\{0,1\}$. Since $\boC_{\beta}$ has $\kappa^+$-c.c.,
so does $\boB_{\beta}$. This implies that
\[\boB=\boB_{\lambda}=\bigcup_{\alpha<\lambda}\boB_{\alpha}\]
because $cf(\lambda)>\kappa$.
Also for each $\beta<\lambda$ one has 
\[|\boC_{\beta}|\leqslant |\beta|^{<\kappa}<\lambda.\] 
Hence $|\boB_{\beta}|\leqslant |\boC_{\beta}|^{\kappa}<\lambda$. 
Next we build an ultrafilter $\E_{\beta}$ on $\boB_{\beta}$ 
for every $\beta\leqslant\lambda$ such that $\omega^{\boB_{\beta}}/\E_{\beta}$
is a proper end-extension of $\omega^{\boB_{\alpha}}/\E_{\alpha}$
when $\alpha<\beta\leqslant\lambda$, and 
\[\omega^{\boB_{\lambda}}/\E_{\lambda}=\bigcup_{\alpha<\lambda}
\omega^{\boB_{\alpha}}/\E_{\alpha}.\]
We first construct retractions
$r_{\alpha,\gamma}:\boB_{\gamma}\mapsto\boB_{\alpha}$ such that
$r_{\alpha,\gamma}=r_{\alpha,\beta}\circ r_{\beta,\gamma}$
for any $\alpha<\beta<\gamma\leqslant\lambda$.
We want also $r_{\alpha,\alpha+1}(a)=0$ for every $a\in I_{\alpha}$,
where $I_{\alpha}$ is the ideal in $\boB_{\alpha+1}$ 
generated by $\{a_{\alpha,n}:n\in\omega\}
\cup\{-\bigvee_{n\in\omega}a_{\alpha,n}\}$.

\medskip

{\bf Claim \ref{arch1}.1.}\quad $I_{\alpha}\cap\boB_{\alpha}=\{0\}$.

\medskip

Proof of Claim \ref{arch1}.1:\quad
Suppose $a\in I_{\alpha}\cap\boB_{\alpha}$. Then there exists
an $m\in\omega$ such that 
\[a\leqslant (\bigvee_{n<m}a_{\alpha,n})\vee
(-\bigvee_{n\in\omega}a_{\alpha,n}).\]
So one has \[a\wedge(-\bigvee_{n<m}a_{\alpha,n})\wedge
(\bigvee_{n\in\omega}a_{\alpha,n})=0.\]
Hence \[a\wedge(\bigwedge_{n<m}-a_{\alpha,n})\wedge a_{\alpha,m}=0.\]
This implies $a=0$ because of the independence of $\{a_{\alpha,n}:n\in\omega\}$
over $\boB_{\alpha}$.\quad $\Box$(Claim \ref{arch1}.1)

\medskip

We construct $r_{\alpha,\beta}$ for any $\alpha<\beta<\delta$ inductively
on $\delta$ when $\delta\leqslant\lambda$. Suppose we have already
$\{r_{\alpha,\beta}:\alpha<\beta<\delta\}$ and want to find 
$r_{\alpha,\delta}$ for every $\alpha<\delta$.
It is trivial when $\delta=0$ or $1$.

Case 1:\quad $\delta$ is a limit ordinal.

Apply Lemma \ref{Shelah2} to obtain retractions 
$p_{\alpha}:\boB_{\delta}\mapsto\boB_{\alpha}$ such that $p_{\alpha}
=r_{\alpha,\beta}\circ p_{\beta}$ for any $\alpha<\beta<\delta$.
Now let $r_{\alpha,\delta}=p_{\alpha}$.

Case 2:\quad $\delta=\beta+1$.

Apply Lemma \ref{Sikorski} to obtain a retraction 
$r:\boB_{\delta}\mapsto\boB_{\beta}$ such that $r(a)=0$ for every 
$a\in I_{\beta}$. Now let $r_{\alpha,\delta}=r_{\alpha,\beta}\circ r$
for every $\alpha<\beta$ and let $r_{\beta,\delta}=r$.
Clearly, the set of retractions $\{r_{\alpha,\delta}:\alpha<\delta\}$
is what we want.

Let $\E_0=\{1\}$ and let $\E_{\alpha}=r^{-1}_{0,\alpha}(\E_{\alpha})$.
It is easy to see that $\E_{\beta}=r^{-1}_{\alpha,\beta}(\E_{\beta})$
for any $\alpha<\beta\leqslant\lambda$.
By Lemma \ref{Koppelberg} $\omega^{\boB_{\beta}}/\E_{\beta}$ 
is an end-extension of $\omega^{\boB_{\alpha}}/\E_{\alpha}$ whenever 
$\alpha<\beta\leqslant\lambda$. Define
$f_{\alpha}\in\omega^{\boB_{\alpha+1}}/\E_{\alpha+1}$ such that
\[f_{\alpha}(0)=-\bigvee_{n\in\omega}a_{\alpha,n}\mbox{ and }
f(n+1)=a_{\alpha,n}\wedge(-\bigvee_{m<n}a_{\alpha,m}).\] Then
$\bigvee_{n\in\omega}r_{\alpha,\alpha+1}(f(n))=0$. Hence
by Lemma \ref{Koppelberg} \[(f_{\alpha})_{\E_{\alpha+1}}\in
\omega^{\boB_{\alpha+1}}/\E_{\alpha+1}\smallsetminus
\omega^{\boB_{\alpha}}/\E_{\alpha},\] {\em i.e.}
$\omega^{\boB_{\alpha+1}}/\E_{\alpha+1}$ is a proper end-extension of
$\omega^{\boB_{\alpha}}/\E_{\alpha}$.
It is also easy to see that 
\[\boul=\bigcup_{\alpha<\lambda} \omega^{\boB_{\alpha}}/\E_{\alpha},\]
where $\boB=\boB_{\lambda}$ and $\E=\E_{\lambda}$,
because $\boB=\bigcup_{\alpha<\lambda}\boB_{\alpha}$
and $cf(\lambda)>\omega$. Hence $|\boul|\geqslant\lambda$.
On the other hand, if $x\in\boul$, then there is a $\beta<\lambda$
such that $x\in\omega^{\boB_{\beta}}/\E_{\beta}$. Hence
\[|\{y\in\boul:y\leqslant x\}|\leqslant|\omega^{\boB_{\beta}}/\E_{\beta}|
\leqslant |\boB_{\beta}|^{\omega}<\lambda.\] This shows that 
$\boul$ is $\lambda$-Archimedean. Now the theorem follows from
Lemma \ref{iso}\quad $\Box$

\bigskip

Next we show a different way to construct $\lambda$-Archimedean ultrapowers.
The construction of a $\lambda$-Archimedean ultrapower in Theorem \ref{arch2}
needs only to assume a measurable cardinal.

\begin{Lemma} \label{embedding}

Suppose $\lambda\geqslant\kappa$ and $\kappa$ is a $\lambda$-supercompact 
cardinal in $V$. Let $\cal U$ be a $\kappa$-complete normal ultrafilter 
on $\kl$ in $V$ and let $j$ be the elementary embedding from $V$ to 
$M=V^{\kl}/{\cal U}$ induced by $\cal U$. Suppose $\boB\subseteq V_{\kappa}$
is a cBa and $j(\boB)\cong\boB*\dot{\boC}$, 
{\em i.e.} $\boB\cong\{j(p):p\in\boB\}$ 
is completely embedded into $j(\boB)$.
Then in $V^{\boB}$ there exists a $\kappa$-complete filter $\D$ on $(\kl)^V$ 
such that ${\cal P}((\kl)^V)/\D\cong\boC$.

\end{Lemma}

\noindent {\bf Proof:}\quad 
Let $G\subseteq\boB$ be a $V$-generic filter and
let $H\subseteq\boC$ be an $M[G]$-generic filter. Let $\hat{j}$ be the 
embedding from $V[G]$ to $M[G][H]$ defined by letting
$\hat{j}(x)=(j(\dot{x}))_{G*H}$. It is well-known that $\hat{j}$ is an 
elementary embedding. Note that $M[G]\subseteq V[G]$ but
$M[G][H]$ is generally not a subclass of $V[G]$. 
For any $A\subseteq (\kl)^V$ in $V[G]$ let 
\[i(A)=[\![j''\lambda\in (j(\dot{A}))_G]\!]_{\boC},\]
where $\dot{A}$ is a $\boB$-name for $A$ and
$j''\lambda=\{j(\alpha):\alpha<\lambda\}\in M^{\lambda}\subseteq M$.

\medskip

{\bf Claim \ref{embedding}.1.}\quad $i$ is a Boolean homomorphism
from ${\cal P}((\kl)^V)$ to $\boC$ in $V[G]$.

\medskip

Proof of Claim \ref{embedding}.1:\quad It suffices to show that
$i$ is well-defined function. The rest follows from the fact that
$j$ is an elementary embedding. 
Suppose $\dot{A}_1$ and $\dot{A}_2$ are two $\boB$-names for 
same set $A$ in $V[G]$. Then there exists a $p\in G$ such that
\[p\Vdash_{\boB}\dot{A}_1=\dot{A}_2.\] By the fact that $j(p)=p$ one has
\[p\Vdash_{\boB*\boC} j(\dot{A}_1)=j(\dot{A}_2)\] or
\[p\Vdash_{\boB}\;\Vdash_{\boC} j(\dot{A}_1)=j(\dot{A}_2).\] 
Since $p\in G$, then 
\[V[G]\models\;\Vdash_{\boC} (j(\dot{A}_1))_G=(j(\dot{A}_2))_G.\] 
Hence in $V[G]$ 
\[[\![j''\lambda\in (j(\dot{A}_1))_G]\!]_{\boC}=
 [\![j''\lambda\in (j(\dot{A}_2))_G]\!]_{\boC}.
\quad\Box(\mbox{(Claim \ref{embedding}.1)}\]

Let's define a filter $\D$ in $V[G]$ by letting 
\[\D=\{A\subseteq(\kl)^V:i(A)=1_{\boC}\}.\]
It is easy to see that $\D\in V[G]$ is a 
$\kappa$-complete filter on $(\kl)^V$.

\medskip

{\bf Claim \ref{embedding}.2.}\quad
${\cal P}((\kl)^V)/\D\cong\boC$ in $V[G]$.

\medskip

Proof of Claim \ref{embedding}.2:\quad
It suffices to show that $i$ is an onto map. Given any $a\in\boC$, let
$a$ be expressed as $[F]_{\cal U}$ for some $F:\kl\mapsto\boB$
in $V$ such that $F$ is not equal to any single $p\in\boB$ modulo $\cal U$.
Let $\dot{A}$ be a $\boB$-name for a subset of $\kl$ in $V$ such that
for each $\sigma\in\kl$ one has 
$[\![\sigma\in\dot{A}]\!]_{\boB}=F(\sigma)$.
Since $j''\lambda=[d]_{\cal U}$, where $d$ is the identity function
from $\kl$ to $\kl$,
\[[\![j''\lambda\in (j(\dot{A}))_G]\!]_{\boC}
=[\![[d]_{\cal U}\in (j(\dot{A}))_G]\!]_{\boC}=
[\langle[\![\sigma\in\dot{A}]\!]_{\boB}:\sigma\in\kl\rangle]_{\cal U}=\]
\[[\langle F(\sigma):\sigma\in\kl\rangle ]_{\cal U}=[F]_{\cal U}=a.\]
Hence $i$ is an onto map from ${\cal P}((\kl)^V)$ to $\boC$.
\quad $\Box$

\begin{Theorem} \label{arch2}

Suppose $\kappa$ is a measurable cardinal such that 
$\theta^{\omega}<2^{\kappa}$ for any $\theta<2^{\kappa}$. Let 
$\boB_{\kappa}$ be the Boolean algebra for adding 
$\kappa$ Cohen reals or adding $\kappa$
random reals. Then in $V^{\boB_{\kappa}}$ there exists an ultrafilter
$\F$ on $\kappa$ such that $\ul$ is $2^{\kappa}$-Archimedean.

\end{Theorem}

\noindent {\bf Proof:}\quad
Obviously, $\kappa$ is $\kappa$-supercompact. Let $j$ be the elementary
embedding induced by a $\kappa$-complete normal ultrafilter $\cal U$
on $\kappa$. For any $\lambda$ let $\boB_{\lambda}$ denote the cBa
for adding $\lambda$ Cohen reals (or random reals). 
Then $j(\boB_{\kappa})=\boB_{j(\kappa)}=\boB_{\kappa}*\dot{\boC}$, where
$\boC$ is isomorphic to the cBa for adding $2^{\kappa}$ Cohen (or random)
reals in $V^{\boB_{\kappa}}$ because 
$|j(\kappa)\smallsetminus\kappa|=2^{\kappa}$.
By Lemma \ref{embedding} there exists a $\kappa$-complete
filter $\D$ on $\kappa$ such that 
\[{\cal P}(\kappa)/\D\cong\boC\] 
in $V^{\boB_{\kappa}}$. 
For any $\lambda$ let $\boB'_{\lambda}$ denote the cBa
for adding $\lambda$ Cohen reals (or random reals) in $V^{\boB_{\kappa}}$. 
Note that $\boC\cong\boB'_{2^{\kappa}}$. 
Since $\boB'_{\lambda}$ has c.c.c. for any
$\lambda$, then $\boB'_{2^{\kappa}}=\bigcup_{\alpha<2^{\kappa}}
\boB'_{\alpha\cdot\omega}$. By the facts that
$|\boB'_{\alpha\cdot\omega}|^{\omega}<2^{\kappa}$
and there exists an independent sequence $\{a_{\alpha,n}:n\in\omega\}$
in $\boB'_{\alpha\cdot\omega+\omega}$ over $\boB'_{\alpha\cdot\omega}$ 
for every $\alpha<2^{\kappa}$, then, by a similar argument in the 
proof of Theorem \ref{arch1}, one can construct an $\F$ such that
$\ul$ is $2^{\kappa}$-Archimedean. \quad $\Box$

\begin{Remark}

(1) Note that $\kappa=2^{\omega}$ in $V^{\boB_{\kappa}}$ and 
it is impossible to have a $2^{\omega}$-Archimedean ultrapower 
by $\omega_1$-saturation. (2) We could have a more general theorem
with a much more general $\boB_{\kappa}$. Restricting $\boB_{\kappa}$ to be
a cBa for adding $\kappa$ Cohen or random reals is just for simplicity
to illustrate the idea.

\end{Remark}

\begin{Theorem}

Suppose $\kappa$ is a supercompact cardinal and $\boB_{\kappa}$ is as
in Theorem \ref{arch2}. Then in $V^{\boB_{\kappa}}$ there exists an 
ultrafilter $\F$ on $\lambda$ for every $\lambda>\kappa$ with 
$\lambda^{<\kappa}=\lambda$ and $\theta^{\omega}<2^{\lambda}$ for every
$\theta<2^{\lambda}$ such that $\omega^{\lambda}/\F$ is 
$2^{\lambda}$-Archimedean.

\end{Theorem}

\noindent {\bf Proof:}\quad Almost identical to the proof of
Theorem \ref{arch2}. \quad $\Box$

\bigskip

Department of Mathematics, College of Charleston,
Charleston, SC 29424

{\em e-mail: jin@@math.cofc.edu}

\medskip

Institute of Mathematics, The Hebrew University,
Jerusalem, Israel

Department of Mathematics, Rutgers University,
New Brunswick, NJ 08903

Department of Mathematics, University of Wisconsin,
Madison, WI 53706

\medskip

{\em Sorting:} The first address is the first author's
and the last three are the second author's.

\end{document}